# Some results on Green's higher Abel-Jacobi map

By Claire Voisin[*]

## 1. Introduction

This paper is devoted to the study of the first higher Abel-Jacobi invariant defined by M. Green in [4] for zero-cycles on a surface. Green's work is a very original attempt to understand, at least over $\mathbb{C}$, the graded pieces of the conjecturally defined filtration on Chow groups

$$\begin{aligned}\mathrm{CH}^p(X)_\mathbb{Q} &= F^0\mathrm{CH}^p(X)_\mathbb{Q} \supset F^1\mathrm{CH}^p(X)_\mathbb{Q} \\ &= \mathrm{CH}^p(X)_\mathbb{Q}^{\mathrm{hom}} \supset \ldots \supset F^{p+1}\mathrm{CH}^p(X)_\mathbb{Q} = 0.\end{aligned}$$

This filtration should satisfy the following properties:

i) First of all it should be stable under correspondences, so that a correspondence $\Gamma \subset X \times Y$ should induce

$$\Gamma_* : F^k\mathrm{CH}^p(X)_\mathbb{Q} \to F^k\mathrm{CH}^{p'}(Y)_\mathbb{Q},$$

where $p' = p + \dim Y - \dim \Gamma$, and

ii) the induced map

$$\mathrm{Gr}^k\Gamma_* : \mathrm{Gr}_F^k\mathrm{CH}^p(X)_\mathbb{Q} \to \mathrm{Gr}_F^k\mathrm{CH}^{p'}(Y)_\mathbb{Q}$$

should vanish when $\Gamma$ is homologous to zero.

A filtration satisfying this last property has been constructed by Saito [11], but it is not shown that the filtration terminates, that is $F^{p+1}\mathrm{CH}^p(X)_\mathbb{Q} = 0$. A definition has also been proposed by J. P. Murre ([8]), under the assumption that a strong Künneth decomposition of the diagonal exists, but it is not proved to satisfy condition ii) above. In fact proving the existence of such a filtration would solve in particular Bloch's conjecture on zero-cycles of surfaces [1].

In any case, the first steps of the filtration are easy to understand, at least for zero-cycles. Namely one should have $F^2\mathrm{CH}_0(X) = \mathrm{CH}_0(X)_{\mathrm{alb}} = \mathrm{Ker}\,\mathrm{alb}$, where $\mathrm{alb} : \mathrm{CH}_0(X)_{\mathrm{hom}} \to \mathrm{Alb}(X)$ is the Albanese map. More generally for the subgroup $\mathrm{CH}^p(X)^{\mathrm{alg}}_{\mathbb{Q}} \subset \mathrm{CH}^p(X)_{\mathbb{Q}}$ of cycles algebraically equivalent to zero, one should have over $\mathbb{C}$, $F^2\mathrm{CH}^p(X)^{\mathrm{alg}}_{\mathbb{Q}} = \mathrm{Ker}\,\Phi^p_X$, where

$$\Phi^p_X : \mathrm{CH}^p(X)^{\mathrm{hom}}_{\mathbb{Q}} \to J^{2p-1}_X = H^{2p-1}(X,\mathbb{C})/F^p H^{2p-1}(X) \oplus H^{2p-1}(X,\mathbb{Z})$$

is the Abel-Jacobi map, defined by Griffiths [5].

---

[*]Partially supported by the project "Algebraic Geometry in Europe" (AGE).

In [4], M. Green suggested constructing directly (over $\mathbb{C}$), from Hodge theoretic considerations, higher Abel-Jacobi maps

$$\psi^p_m : F^m\mathrm{CH}^p(X) \to J^p_m(X),$$

so that $F^{m+1}\mathrm{CH}^p(X) = \mathrm{Ker}\,\psi^p_m$. Hence one should have an induced injective map

$$\psi^p_m : \mathrm{Gr}^m\mathrm{CH}^p(X) \to J^p_m(X).$$

In the case of zero cycles on a surface, he proposed an explicit construction of

(1.1) $$\psi^2_2 : \mathrm{Gr}^2\mathrm{CH}^2(S) = \mathrm{CH}_0(S)_{\mathrm{alb}} \to J^2_2(S)$$

that we will review below. The purpose of this paper is to answer some questions raised in [4], concerning the behaviour of $\psi^2_2$. To simplify the notation, we will assume throughout that $S$ is regular, but this assumption does not play any role in the arguments. Our first result is the following, which answers negatively conjecture 3.4 of [4]:

THEOREM 1. *The higher Abel-Jacobi map $\psi^2_2$ is not, in general, injective.*

The noninjectivity is proved here for an explicit example but the argument should allow us to prove, as we will explain, that $\psi^2_2$ is never injective for surfaces with $\mathrm{CH}_0(S)_{\mathrm{alb}} \neq 0$.

Our second result solves, in particular, conjecture 3.6 of [4]:

THEOREM 2. *The map $\psi^2_2$ is nontrivial modulo torsion (and has an infinite dimensional image), when $h^{2,0}(S) \neq 0$.*

As an intermediate step, we explain how Mumford's pull-back $Z^*(\omega)$ of a holomorphic two-form $\omega$ of a surface $S$ on a variety $W$ parametrizing 0-cycles $(Z_w)_{w \in W}$ of $S$ (cf. [7]) can be computed when one has a family $\mathcal{C} \to W$ of curves of $S$ parametrized by $W$, such that for each $w \in W$, the 0-cycle $Z_w$ is supported on $C_w$. There are then two associated Abel-Jacobi invariants $e_{C_w,S}$ and $f_{Z_w,C_w}$ to be defined below, which play a key role in the construction of $\psi^2_2(Z_w)$, and we show that $Z^*(\omega)$ can be computed from the wedge product $de \wedge df$. We then use this result to show that in fact $Z^*\omega$ depends only on the map $\psi^2_2 \circ Z_* : W \to J^2_2(S)$.

Thus our results show that the first new higher Abel-Jacobi map defined by Mark Green is not strong enough to capture the whole of $\mathrm{CH}_0(S)_{\mathrm{alb}}$ as it should conjecturally do, but that it is strong enough to determine Mumford's invariants, which were used to show that $\mathrm{CH}_0(S)_{\mathrm{alb}}$ is infinite dimensional, when $h^{2,0}(S) \neq 0$. The question of whether it is possible to refine it so as to get the desired injectivity of 1.1 is still open.

The paper is organized as follows: The result above concerning the pullback of holomorphic two forms (Proposition 2) provides the contents of Section 3. Theorem 1 is proved in Section 2, and Theorem 2 is proved in Section 4.

We conclude this introduction with a brief description of $\psi_2^2$, which will serve also as an introduction for the notation used throughout the paper.

Let $S$ be a regular surface, and let $C$ be a smooth (not necessarily connected) curve; let $\psi : C \to S$ be a morphism generically one-to-one on its image. We can find an immersion $\phi : C \hookrightarrow \tilde{S}$, and a birational morphism $\tau : \tilde{S} \to S$ such that $\psi = \tau \circ \phi$.

Now let $Z$ be a 0-cycle of $C$, of degree 0 on each component of $C$. We construct two Abel-Jacobi invariants $e_{C,S}$ and $f_{Z,C}$ as follows:

The mixed Hodge structure on $H^2(\tilde{S}, C)$ is given by the Hodge filtration $F^{\cdot}$ on $H^2(\tilde{S}, C)$, which fits in the exact sequence

$$(1.2) \qquad 0 \to H^1(C, \mathbb{Z}) \to H^2(\tilde{S}, C, \mathbb{Z}) \to \mathrm{Ker}(H^2(\tilde{S}, \mathbb{Z}) \to H^2(C, \mathbb{Z})) \to 0.$$

The filtration $F^{\cdot}$ restricts to the Hodge filtration on $H^1(C)$ and projects to the Hodge filtration on $\mathrm{Ker}(H^2(\tilde{S}, \mathbb{Z}) \to H^2(C, \mathbb{Z}))$.

Define $H^2(S, \mathbb{Z})_{\mathrm{tr}}$ as the quotient $H^2(S, \mathbb{Z})/NS(S)$. Then its dual $H^2(S, \mathbb{Z})_{\mathrm{tr}}\check{\phantom{t}}$ is the orthogonal of $NS(S)$ in $H^2(S, \mathbb{Z})$. There is an inclusion of Hodge structures

$$\tau^* : H^2(S, \mathbb{Z})_{\mathrm{tr}}\check{\phantom{t}} \hookrightarrow \mathrm{Ker}(H^2(\tilde{S}, \mathbb{Z}) \to H^2(C, \mathbb{Z})).$$

Restricting the extension (1.2) to $H^2(S, \mathbb{Z})_{\mathrm{tr}}\check{\phantom{t}}$, we get an exact sequence of mixed Hodge structures

$$0 \to H^1(C, \mathbb{Z}) \to H^2(\tilde{S}, C, \mathbb{Z})_{\mathrm{tr}} \to H^2(S, \mathbb{Z})_{\mathrm{tr}}\check{\phantom{t}} \to 0.$$

The extension class of this exact sequence is an element $e_{C,S}$ of the complex torus (cf. [3]),

$$\begin{aligned} J(C \times S)_{\mathrm{tr}} &:= H^1(C, \mathbb{C}) \otimes H^2(S, \mathbb{C})_{\mathrm{tr}} / [F^2(H^1(C) \otimes H^2(S)_{\mathrm{tr}}) \\ &\quad \oplus H^1(C, \mathbb{Z}) \otimes H^2(S, \mathbb{Z})_{\mathrm{tr}}]. \end{aligned}$$

It is not difficult to prove that it can be also computed as the natural projection of the Abel-Jacobi invariant of the one-cycle obtained from the graph of $\psi$ (which is a one-cycle of $C \times S$) by adding vertical and horizontal one-cycles of $C \times S$ in order to get a homologically trivial one-cycle.

It is well-known that the inclusion

$$H^1(C, \mathbb{R}) \otimes H^2(S, \mathbb{R})_{\mathrm{tr}} \subset H^1(C, \mathbb{C}) \otimes H^2(S, \mathbb{C})_{\mathrm{tr}}$$

induces an isomorphism

$$H^1(C, \mathbb{R}) \otimes H^2(S, \mathbb{R})_{\mathrm{tr}} \cong H^1(C, \mathbb{C}) \otimes H^2(S, \mathbb{C})_{\mathrm{tr}} / F^2(H^1(C) \otimes H^2(S)_{\mathrm{tr}}),$$

and this allows us to identify $J(C \times S)_{\text{tr}}$ to the real torus

$$H^1(C,\mathbb{Z}) \otimes_\mathbb{Z} H^2(S,\mathbb{Z})_{\text{tr}} \otimes_\mathbb{Z} \mathbb{R}/\mathbb{Z}.$$

We will view $e_{C,S}$ as an element of this real torus.

Now the zero-cycle $Z$ has an Abel-Jacobi invariant (Albanese image)

$$f_{Z,C} \in J(C) \cong H^1(C,\mathbb{C})/[F^1 H^1(C) \oplus H^1(C,\mathbb{Z})].$$

Again, the inclusion $H^1(C,\mathbb{R}) \subset H^1(C,\mathbb{C})$ induces an isomorphism $H^1(C,\mathbb{R}) \cong H^1(C,\mathbb{C})/F^1 H^1(C)$, which provides the identification

$$J(C) \cong H^1(C,\mathbb{Z}) \otimes_\mathbb{Z} \mathbb{R}/\mathbb{Z}.$$

We will view $f_{Z,C}$ as an element of the real torus on the right.

The pairing

$$H^1(C,\mathbb{Z}) \otimes_\mathbb{Z} H^1(C,\mathbb{Z}) \to \mathbb{Z}$$

allows us then to contract $e_{C,S}$ and $f_{Z,C}$ to an element

$$e_{C,S} \cdot f_{Z,C} \in \mathbb{R}/\mathbb{Z} \otimes_\mathbb{Z} \mathbb{R}/\mathbb{Z} \otimes_\mathbb{Z} H^2(S,\mathbb{Z})_{\text{tr}}.$$

Defining now $U_2^2 \subset \mathbb{R}/\mathbb{Z} \otimes_\mathbb{Z} \mathbb{R}/\mathbb{Z} \otimes_\mathbb{Z} H^2(S,\mathbb{Z})_{\text{tr}}$ as the group generated by the elements $e_{C,S} \cdot f_{Z,C}$ defined above, for the triples $(C, Z, \psi)$ such that $\psi_*(Z) = 0$ as a zero-cycle of $S$, it is clear that the projection

$$\overline{e_{C,S} \cdot f_{Z,C}} \in J_2^2(S) := \mathbb{R}/\mathbb{Z} \otimes_\mathbb{Z} \mathbb{R}/\mathbb{Z} \otimes_\mathbb{Z} H^2(S,\mathbb{Z})_{\text{tr}}/U_2^2$$

depends only on the zero-cycle $\psi_*(Z)$. The resulting map

$$\psi_2^2 : Z_0(S)_{\text{hom}} \to J_2^2(S)$$

is then easily seen to factor through rational equivalence, so that $\psi_2^2$ is actually defined on $\text{CH}_0^0(S)$.

*Acknowledgements.* I would like to thank P. Griffiths for his careful reading of the first version of this paper, and for the improvements he suggested.

## 2. The noninjectivity of $\psi_2^2$

In this section we construct a counterexample to the conjectured injectivity of the map

$$\psi_2^2 : \text{CH}_0(S)_{\text{alb}} \to J_2^2(S).$$

The counterexample is based on a refinement of the following argument.

First of all, if $\Gamma \subset C \times S$ is a correspondence homologous to zero, with Abel-Jacobi invariant

$$e_\Gamma \in J(C \times S)_{\text{tr}},$$

we show that
$$\psi_2^2 \circ \Gamma_* : CH_0(C)_{\text{hom}} \to J_2^2(S)$$
is given by
$$\psi_2^2 \circ \Gamma_*(z) = e_\Gamma \cdot f_z \bmod U_2^2,$$
where $f_z = \text{alb}(z) \in J(C) = H^1(C, \mathbb{Z}) \otimes \mathbb{R}/\mathbb{Z}$. Now we view $e_\Gamma$ as an element of $\text{Hom}(H^2(S, \mathbb{Z})_{\text{tr}}\check{}, H^1(C, \mathbb{Z}) \otimes \mathbb{R}/\mathbb{Z})$ and we note that if its image is contained in a proper real subtorus $T$ of $H^1(C, \mathbb{Z}) \otimes \mathbb{R}/\mathbb{Z}$, there is a nontrivial real subtorus $T^\perp$ of $H^1(C, \mathbb{Z}) \otimes \mathbb{R}/\mathbb{Z}$ such that, if $f_z \in T^\perp$,
$$e_\Gamma \cdot f_z = 0 \text{ in } H^2(S, \mathbb{Z})_{\text{tr}} \otimes \mathbb{R}/\mathbb{Z} \otimes \mathbb{R}/\mathbb{Z}.$$
Then the injectivity of $\psi_2^2$ would imply that $T^\perp \subset \text{Ker}\,\Gamma_*$, and if $J(C)$ is simple, this would imply that $\Gamma_* = 0$, and then $e_\Gamma$ would be a torsion point in $J(C \times S)_{\text{tr}}$. So it suffices to find $C, \Gamma$ as above with $e_\Gamma$ not of torsion (or $\Gamma_* \neq 0$), $J(C)$ simple and $\text{Im}\,e_\Gamma$ contained in a proper real subtorus $T$ of $H^1(C, \mathbb{Z}) \otimes \mathbb{R}/\mathbb{Z}$ to contradict the injectivity of $\psi_2^2$.

We start with the simple Lemma 1 below which allows us to extend slightly the definition of $\psi_2^2$. Let $S$ be a regular surface, $C$ be a smooth curve and $\Gamma \in CH_1(C \times S)$ be a one-cycle; the homology class of $\Gamma$ lies in $H_2(C) \oplus H_2(S)_{\text{alg}}$, so that adding to $\Gamma$ vertical and horizontal cycles we can get a cycle $\Gamma'$ homologous to zero: Then the induced morphisms
$$\Gamma_* : CH_0^0(C) \to CH_0^0(S), \quad \Gamma'_* : CH_0^0(C) \to CH_0^0(S)$$
coincide, and the Abel-Jacobi image of $\Gamma'$ in $J(C \times S)_{\text{tr}}$ (see Section 1) does not depend on the choice of $\Gamma'$. We will denote it by $e_\Gamma$. As in Section 1, we can view $e_\Gamma$ as an element of the real torus
$$H^1(C, \mathbb{Z}) \otimes_\mathbb{Z} H^2(S, \mathbb{Z})_{\text{tr}} \otimes_\mathbb{Z} \mathbb{R}/\mathbb{Z}.$$
By contraction and use of the intersection pairing on $H^1(C, \mathbb{Z})$, $e_\Gamma$ gives a map
$$[e_\Gamma] : JC \cong H^1(C, \mathbb{Z}) \otimes_\mathbb{Z} \mathbb{R}/\mathbb{Z} \to \mathbb{R}/\mathbb{Z} \otimes_\mathbb{Z} \mathbb{R}/\mathbb{Z} \otimes_\mathbb{Z} H^2(S, \mathbb{Z})_{\text{tr}}.$$
We have now:

LEMMA 1. *For $z \in JC$, $\psi_2^2(\Gamma_*(z)) \in J_2^2(S)$ is equal to the projection of $[e_\Gamma](z)$ modulo $U_2^2(S)$, using the definition*
$$J_2^2(S) := \mathbb{R}/\mathbb{Z} \otimes_\mathbb{Z} \mathbb{R}/\mathbb{Z} \otimes_\mathbb{Z} H^2(S, \mathbb{Z})_{\text{tr}}/U_2^2$$
*of Section* 1.

*Proof.* This is true by definition if $\Gamma$ is the graph $\Gamma_\phi$ of a morphism $\phi$ from $C$ to $S$, generically one-to-one on its image. Now let $C_1 \xrightarrow{\phi} S$ be the desingularization of the inclusion of $\mathrm{pr}_2(\mathrm{Supp}\,\Gamma)$ in $S$. Then $\Gamma$ lifts to a one-cycle $\Gamma_1 \in C \times C_1$, such that

$$\Gamma_* = \phi_* \circ \Gamma_{1*} : JC \to \mathrm{CH}_0^0(S).$$

We have then

$$\begin{aligned}\psi_2^2(\Gamma_*(z)) &= \psi_2^2(\phi_*(\Gamma_{1*}(z))) \\ &= \text{projection of } [e_{\Gamma_\phi}](\Gamma_{1*}(z)) \text{ in } J_2^2(S).\end{aligned}$$

Now it suffices to prove that

$$[e_\Gamma] = [e_{\Gamma_\phi}] \circ \Gamma_{1*} : J(C) \to \mathbb{R}/\mathbb{Z} \otimes_\mathbb{Z} \mathbb{R}/\mathbb{Z} \otimes_\mathbb{Z} H^2(S,\mathbb{Z})_{\mathrm{tr}}.$$

But $\Gamma_1$ induces naturally a correspondence $\tilde{\Gamma}_1$ between $C \times S$ and $C_1 \times S$; hence a morphism

$$\tilde{\Gamma}_1^* : \mathrm{CH}_1(C_1 \times S) \to \mathrm{CH}_1(C \times S),$$

such that $\Gamma \equiv_{\mathrm{rat}} \tilde{\Gamma}_1^*(\Gamma_\phi)$. It follows that $e_\Gamma = \tilde{\Gamma}_1^*(e_{\Gamma_\phi})$ in $J(C \times S)_{\mathrm{tr}}$, where $\tilde{\Gamma}_1^*$ also denotes the induced morphism between the intermediate jacobians $J(C_1 \times S)_{\mathrm{tr}}$ and $J(C \times S)_{\mathrm{tr}}$.

Let $\Gamma_{1\mathbb{Z}}^* : H^1(C_1, \mathbb{Z}) \to H^1(C, \mathbb{Z})$ be the morphism of Hodge structures induced by the cohomology class of $\Gamma_1$ in $C \times C_1$; then the morphism $\tilde{\Gamma}_1^*$ is induced by the morphism of Hodge structures

$$\Gamma_{1\mathbb{Z}}^* \otimes \mathrm{Id} : H^1(C_1, \mathbb{Z}) \otimes H^2(S, \mathbb{Z})_{\mathrm{tr}} \to H^1(C, \mathbb{Z}) \otimes H^2(S, \mathbb{Z})_{\mathrm{tr}},$$

and it follows that we have a commutative diagram

$$\begin{array}{lllll}\Gamma_{1\mathbb{R}}^* \otimes \mathrm{Id} & : H^1(C_1, \mathbb{R}) \otimes H^2(S, \mathbb{R})_{\mathrm{tr}} & \to & H^1(C, \mathbb{R}) \otimes H^2(S, \mathbb{R})_{\mathrm{tr}} \\ & \downarrow & & \downarrow \\ \tilde{\Gamma}_{1\mathbb{C}}^* \bmod F^2 & : H^1(C_1, \mathbb{C}) \otimes H^2(S, \mathbb{C})_{\mathrm{tr}}/F^2 & \to & H^1(C, \mathbb{C}) \otimes H^2(S, \mathbb{C})_{\mathrm{tr}}/F^2,\end{array}$$

where the vertical arrows are the identifications already used between real cohomology and complex cohomology mod $F^2$, and the last horizontal map induces

$$\tilde{\Gamma}_1^* : J(C_1 \times S)_{\mathrm{tr}} \to J(C \times S)_{\mathrm{tr}}$$

by passing to the quotient modulo integral cohomology. This means that, viewed as elements of $J(C_1) \otimes_\mathbb{Z} H^2(S, \mathbb{Z})_{\mathrm{tr}}$ and $J(C) \otimes_\mathbb{Z} H^2(S, \mathbb{Z})_{\mathrm{tr}}$ respectively, $e_{\Gamma_\phi}$ and $e_\Gamma$ satisfy the relation

$$e_\Gamma = \Gamma_1^* \otimes \mathrm{Id}(e_{\Gamma_\phi}).$$

Now it suffices to note that the contraction maps

$$\begin{array}{ll}\langle,\rangle_{C_1} : & J(C_1) \otimes_\mathbb{Z} J(C_1) \to \mathbb{R}/\mathbb{Z} \otimes_\mathbb{Z} \mathbb{R}/\mathbb{Z}, \\ \langle,\rangle_C : & J(C) \otimes_\mathbb{Z} J(C) \to \mathbb{R}/\mathbb{Z} \otimes_\mathbb{Z} \mathbb{R}/\mathbb{Z}\end{array}$$

satisfy the relation
$$\langle \Gamma_{1*}(z), w \rangle_{C_1} = \langle z, \Gamma_1^*(w) \rangle_C, \ z \in J(C), \ w \in J(C_1),$$
to get
$$\begin{aligned}[e_{\Gamma_\phi}](\Gamma_{1*}(z)) &= \langle \Gamma_{1*}(z), e_{\Gamma_\phi} \rangle_{C_1} \\ &= \langle z, \Gamma_1^* \otimes \mathrm{Id}(e_{\Gamma_\phi}) \rangle_C = \langle z, e_\Gamma \rangle_C = [e_\Gamma](z),\end{aligned}$$
as desired. $\square$

The following Lemma 2 is quite standard (cf. [10]); let $\Gamma \in \mathrm{CH}_1(C \times S)$ be a correspondence, and let
$$\Gamma_* : J(C) \to \mathrm{CH}_0^0(S)$$
be the induced morphism; we have:

LEMMA 2. $\mathrm{Ker}\,\Gamma_*$ *is a countable union of translates of an abelian subvariety of* $J(C)$.

*Proof.* $\mathrm{Ker}\,\Gamma_*$ is a subgroup of $J(C)$, and is a countable union of algebraic subsets of $J(C)$. The union of the irreducible algebraic subsets of $J(C)$ passing through 0 and contained in $\mathrm{Ker}\,\Gamma_*$ is stable under difference which implies that it can be written as an increasing union of irreducible algebraic subsets of $J(C)$. So it must be in fact an algebraic subset of $J(C)$, stable under difference, that is an abelian subvariety of $J(C)$. Hence the result. $\square$

Now assume some real subtorus $T$ of $J(C) = H^1(C, \mathbb{R})/H^1(C, \mathbb{Z})$ is contained in $\mathrm{Ker}\,\Gamma_*$; then if $A \subset J(C)$ is the maximal abelian subvariety contained in $\mathrm{Ker}\,\Gamma_*$, so that by Lemma 2, $\mathrm{Ker}\,\Gamma_* = \bigcup_{m \in \mathbb{Z}} A + t_m$ for some $t_m \in J(C)$, then
$$T = \bigcup_{m \in \mathbb{Z}} T \cap (A + t_m).$$
It follows that some $T \cap (A + t_m)$ must contain an open set of $T$, and this implies easily that in fact $T$ is contained in $A$. So we have proved:

LEMMA 3. *Let $T$ be a real subtorus of $J(C)$ contained in $\mathrm{Ker}\,\Gamma_*$; then there is an abelian subvariety $A$ of $J(C)$ such that $T \subset A \subset \mathrm{Ker}\,\Gamma_*$. In particular, if $T$ is nontrivial and $T \subset B$ where $B$ is a simple abelian subvariety of $J(C)$ (i.e. there is no proper nontrivial abelian subvariety of $B$), then $B \subset \mathrm{Ker}\,\Gamma_*$.*

We want to apply these observations to show the noninjectivity of the higher Abel-Jacobi map $\psi_2^2 : \mathrm{CH}_0^0(S) \to J_2^2(S)$. Let $C$ be a curve, and $\Gamma \in \mathrm{CH}_1(C \times S)$ be a correspondence. Let $e_\Gamma \in J(C) \otimes_{\mathbb{Z}} H^2(S, \mathbb{Z})_{\mathrm{tr}}$ be the corresponding Abel-Jacobi invariant. We can view $e_\Gamma$ as an element
$$[e_\Gamma]^* \in \mathrm{Hom}(H^2(S, \mathbb{Z})_{\mathrm{tr}}^{\vee}, J(C)).$$

Assume there is a proper real subtorus $T$ of $J(C)$ containing $\text{Im}\,[e_\Gamma]^*$; i.e. there is a proper sublattice $T_\mathbb{Z}$ of $H^1(C,\mathbb{Z})$ such that $\text{Im}\,[e_\Gamma]^* \subset T_\mathbb{Z} \otimes_\mathbb{Z} \mathbb{R}/\mathbb{Z}$. Then by definition of $[e_\Gamma]$

$$T^\perp := T_\mathbb{Z}^\perp \otimes_\mathbb{Z} \mathbb{R}/\mathbb{Z} \subset \text{Ker}\,[e_\Gamma].$$

Similarly, if $B \xhookrightarrow{i_B} J(C)$ is an abelian subvariety, and $\check{B}$ is the corresponding quotient of $J(C)$, let $J(C \times S)^B_{\text{tr}}$ be the induced quotient of $J(C \times S)_{\text{tr}}$; that is, writing $\check{B} = \check{B}_\mathbb{C}/\check{B}^{1,0} \oplus \check{B}_\mathbb{Z}$, then

$$J(C \times S)^B_{\text{tr}} := \check{B}_\mathbb{C} \otimes H^2(S,\mathbb{C})_{\text{tr}}/F^2(\check{B}_\mathbb{C} \otimes H^2(S,\mathbb{C})_{\text{tr}}) \oplus \check{B}_\mathbb{Z} \otimes H^2(S,\mathbb{Z})_{\text{tr}}.$$

Let

$$[e_\Gamma]^*_B \in \text{Hom}(H^2(S,\mathbb{Z})_{\text{tr}}, \check{B})$$

be the composition of $[e_\Gamma]^*$ with the projection $J(C) \to \check{B}$. Let $e^B_\Gamma \in J(C \times S)^B_{\text{tr}}$ be the projection of $e_\Gamma$. Note that $[e_\Gamma]^*_B$ is simply $e^B_\Gamma$ viewed as an element of $\text{Hom}(H^2(S,\mathbb{Z})_{\text{tr}}, \check{B})$ using the real representations of the (intermediate) jacobians

$$J(C \times S)^B_{\text{tr}} \cong \text{Hom}_\mathbb{Z}(H^2(S,\mathbb{Z})_{\text{tr}}, \check{B}_\mathbb{Z} \otimes_\mathbb{Z} \mathbb{R}/\mathbb{Z}).$$

If $\text{Im}\,[e_\Gamma]^*_B$ is contained in a proper real subtorus $T$ of $\check{B}$, the orthogonal torus $T^\perp \subset B$ is contained in $\text{Ker}\,[e_\Gamma]_{|B}$.

In this situation, assume now that $\psi_2^2$ is injective and that $B$ is simple: then by Lemma 1, one finds that $\Gamma_*$ vanishes on $T^\perp \subset B$, and by Lemma 3, one concludes that $\Gamma_*$ vanishes on $B$. Now this implies:

PROPOSITION 1. *Under the above assumptions, the projection $e^B_\Gamma$ of $e_\Gamma$ in $J(C \times S)^B_{\text{tr}}$ is in fact a torsion point.*

This follows from $\Gamma_{*|B} = 0$ and from the following lemma (cf. [4], [2]) applied to the correspondence $\Gamma \circ \pi_B$ where $\pi_B$ is a multiple of a projector from $J(C)$ to $B$:

LEMMA 4. *Let $\Gamma \in \text{CH}_1(C \times S)$ be a correspondence such that the corresponding map $\Gamma_* : J(C) \to \text{CH}_0^0(S)$ is zero; then the Abel-Jacobi invariant $e_\Gamma$ is a torsion point of $J(C \times S)_{\text{tr}}$.*

In order to contradict the injectivity of $\psi_2^2$ it suffices then to find a smooth curve $C$, a simple abelian subvariety $B$ of $J(C)$ and a correspondence $\Gamma \in \text{CH}_1(C \times S)$ satisfying the following properties:

– The projection $e^B_\Gamma$ of the Abel-Jacobi invariant $e_\Gamma \in J(C \times S)_{\text{tr}}$ in $J(C \times S)^B_{\text{tr}}$ is not a torsion point.

– The image of the map

$$[e_\Gamma]^*_B : H^2(S,\mathbb{Z})_{\text{tr}} \to \check{B}$$

is contained in a proper real subtorus of $\check{B}$.

To get an explicit example, we use a construction due to Paranjape ([9]). Consider a $K3$ surface $S$ which is the desingularization of a general double cover of $\mathbb{P}^2$ branched along the union of six lines. Then $\operatorname{rk} NS(S) = 16$, hence $b_2(S)_{\operatorname{tr}} = 6$. Paranjape constructs a genus 5 curve $C$, which is a ramified cover of an elliptic curve $E$, with an automorphism $j$ of order 4, acting on $B := (\operatorname{Ker} Nm : J(C) \to J(E))^0$, a four dimensional abelian variety. The $K3$ surface $S$ is birational to a quotient of $C \times C$ by a finite group. Let $r : C \times C \to S$ be the quotient (rational) map; for generic $c \in C$, $r$ is everywhere defined along $c \times C$ and we get a family of one-cycles of $C \times S$ parametrized by $C$,

$$c \in C_{\operatorname{gen}} \mapsto \Gamma_c := \text{ graph of } r_{|c \times C} \subset C \times S.$$

This family induces an Abel-Jacobi map

$$\Gamma_* : J(C) \to J(C \times S).$$

Using the projection

$$J(C \times S) \to J(C \times S)_{\operatorname{tr}} \to J(C \times S)_{\operatorname{tr}}^B$$

and restricting the map $\Gamma_*$ to $B \subset J(C)$, we get a morphism (of complex tori)

$$\Gamma_*^B : B \to J(C \times S)_{\operatorname{tr}}^B.$$

This morphism corresponds to a morphism of Hodge structures

$$\phi_\Gamma : B_\mathbb{Z} \to \check{B}_\mathbb{Z} \otimes H^2(S, \mathbb{Z})_{\operatorname{tr}},$$

where $B = B_\mathbb{C}/F^1 B_\mathbb{C} \oplus B_\mathbb{Z}$, $B_\mathbb{Z} = \check{\check{B}}_\mathbb{Z}$. One verifies easily that the corresponding morphism of Hodge structures

$$\psi_\Gamma : H^2(S, \mathbb{Z})_{\operatorname{tr}}\check{} \to \check{B}_\mathbb{Z} \otimes \check{B}_\mathbb{Z}$$

is the composite of the pull-back map

$$r^* : H^2(S, \mathbb{Z})_{\operatorname{tr}}\check{} \to H^1(C, \mathbb{Z}) \otimes H^1(C, \mathbb{Z})$$

and of the projection map

$$H^1(C, \mathbb{Z}) \otimes H^1(C, \mathbb{Z}) \to \check{B}_\mathbb{Z} \otimes \check{B}_\mathbb{Z}.$$

Paranjape ([9]) shows that $\psi_\Gamma$ is injective. It follows that $\phi_\Gamma$ is nonzero.

Now let $u \in B_\mathbb{Z}$ be such that $\phi_\Gamma(u) \neq 0$. There are at most countably many points $u_i$ in the real torus $(\mathbb{R}/\mathbb{Z}) \cdot u$ such that $\Gamma_*^B(u_i)$ is of torsion in $J(C \times S)_{\operatorname{tr}}^B$. Let $\alpha \in \mathbb{R}/\mathbb{Z}$ be such that $\Gamma_*^B(\alpha \cdot u)$ is not of torsion; now view $\phi_\Gamma(u)$ as an element $[\phi_\Gamma(u)]$ of $\operatorname{Hom}(H^2(S, \mathbb{Z})_{\operatorname{tr}}\check{}, \check{B}_\mathbb{Z})$. Since $b_{2,tr}(S) = 6$ and $\operatorname{rk} \check{B}_\mathbb{Z} = 8$, the image of $[\phi_\Gamma(u)]$ is contained in a proper sublattice of $\check{B}_\mathbb{Z}$. It follows that $[\phi_\Gamma(u)] \otimes \alpha \in \operatorname{Hom}(H^2(S, \mathbb{Z})_{\operatorname{tr}}\check{}, \check{B}_\mathbb{Z} \otimes \mathbb{R}/\mathbb{Z})$ has its image contained in a proper subtorus of $\check{B}_\mathbb{Z} \otimes \mathbb{R}/\mathbb{Z}$.

Since $\Gamma_*^B$ is induced by the Abel-Jacobi map, there is a one-cycle $\Gamma_{u \cdot \alpha}$ in $C \times S$ such that $e_{\Gamma_{u \cdot \alpha}}^B = \Gamma_*^B(\alpha \cdot u)$. Consider now the corresponding element $[e_{\Gamma_{u\cdot\alpha}}]^*_B$ of $\mathrm{Hom}(H^2(S,\mathbb{Z})_{\mathrm{tr}}, \check{B}_\mathbb{Z} \otimes \mathbb{R}/\mathbb{Z})$. Since $\phi_\Gamma$ is a morphism of Hodge structures, we have a commutative diagram

$$\begin{array}{ccc} \phi_\Gamma \otimes \mathbb{R}/\mathbb{Z}: & B_\mathbb{Z} \otimes \mathbb{R}/\mathbb{Z} & \to & \mathrm{Hom}(H^2(S,\mathbb{Z})_{\mathrm{tr}}, \check{B}_\mathbb{Z} \otimes \mathbb{R}/\mathbb{Z}) \\ & \downarrow & & \downarrow \\ \Gamma^B_*: & B & \to & J(C \times S)^B_{\mathrm{tr}} \end{array}$$

where the vertical maps are the identifications already used above. It follows that $[e_{\Gamma_{u\cdot\alpha}}]^*_B$ is equal to $[\phi_\Gamma(u)] \otimes \alpha$, hence has its image contained in a proper subtorus of $\check{B}_\mathbb{Z} \otimes \mathbb{R}/\mathbb{Z}$.

To conclude that this is the desired counterexample, it suffices to note that for general $S$, $B$ is a simple abelian variety. This follows from the fact that $(B, j)$ determines the Hodge structure on $H^2(S)_{\mathrm{tr}}$ (cf. [9]), which implies that $B$ has at least four moduli. Then a dimension count shows that the moduli space of nonsimple abelian varieties $A$ of dimension 4 admitting an automorphism of order 4, acting on $H^{1,0}(A)$ with two eigenvalues equal to $i$ and two eigenvalues equal to $-i$, as is the case in Paranjape's family, is of dimension strictly less than 4. □

The counterexample given here is quite special, but it seems from the line of the argument that the noninjectivity of $\psi_2^2$ for a surface with infinite dimensional $\mathrm{CH}_0$ is a general fact; indeed take any such surface $S$ (regular for simplicity) and choose a finite sufficiently ample and generic morphism $\phi: S \to \mathbb{P}^2$. Let $C$ be a sufficiently general and ample curve in $\mathbb{P}^2$ such that $\tilde{C} = \phi^{-1}(C)$ is smooth, $J(C)$ is simple, and $j_*: B \to \mathrm{CH}_0(S)$ has an at most countable kernel, where $j$ is the inclusion of $\tilde{C}$ in $S$ and $B := (\mathrm{Ker}\, Nm: J(\tilde{C}) \to J(C))_0$. Now choose a dimension-1 real subtorus $T$ of $\phi^*(J(C))$ and let $T^\perp \subset J(\tilde{C})$ be its orthogonal. Consider a general small deformation $\tilde{C}_t$ of $\tilde{C}$. The associated element $e_{\tilde{C}_t,S}$ of $J(\tilde{C}_t \times S)_{\mathrm{tr}}$ varies holomorphically with $t$ and the corresponding element $[e_{\tilde{C}_t,S}]^* \in \mathrm{Hom}(H^2(S,\mathbb{Z})_{\mathrm{tr}}, H^1(\tilde{C}_t,\mathbb{Z}) \otimes \mathbb{R}/\mathbb{Z}) \cong \mathrm{Hom}(H^2(S,\mathbb{Z})_{\mathrm{tr}}, H^1(\tilde{C},\mathbb{Z}) \otimes \mathbb{R}/\mathbb{Z})$ varies in a real analytic way. By construction, we have $\mathrm{Im}\,[e_{\tilde{C}_0,S}]^* \subset T^\perp$, and the locus where $\mathrm{Im}\,[e_{\tilde{C}_t,S}]^*$ remains contained in $T^\perp$ is defined by $b_2(S)_{\mathrm{tr}}$ real analytic equations. Now, the arguments developed above show that if $\psi_2^2$ is injective, for $t$ in this locus, there is an abelian subvariety $A_t$ of $J(\tilde{C}_t)$ such that

$$T \subset A_t \subset \mathrm{Ker}\, j_{t*}.$$

The simplicity of $J(C)$ and the fact that $\phi^*(J(C))$ is the maximal abelian subvariety of $J(\tilde{C}_0)$ contained in $\mathrm{Ker}\, j_{0*}$ imply now that on a connected positive dimensional component of this locus containing 0, $A_t \subset J(\tilde{C}_t)$ is a deformation of $\phi^*(J(C)) \subset J(\tilde{C}_0)$.

A contradiction would follow by proving the following facts:

- The small deformations $\tilde{C}_t$ of $\tilde{C} = \tilde{C}_0$ such that $J(\tilde{C}_t)$ contains a deformation of $\phi^*(J(C))$ are the curves of the form $\phi_t^{-1}(C_t)$ where $\phi_t$ is a deformation of $\phi$ and $C_t$ is a deformation of $C$. In particular they form a sublocus of the family of deformations of $\tilde{C}_t$ of arbitrarily large codimension.

- The locus where $\mathrm{Im}\,[e_{\tilde{C}_t,S}]^*$ remains contained in $T^\perp$ is actually of real codimension less or equal to $b_2(S)_{\mathrm{tr}}$. (This is not clear since the equations are only real analytic, and not holomorphic, but this could be proved by an infinitesimal study: it would suffice to show that the equations have independent differentials at 0.)

## 3. A formula for the pull-back of holomorphic two-forms

Let $S$ be a regular surface. Let $W$ be a complex ball parametrizing the following data:

$\mathcal{C}$ is a smooth complex variety, $\pi : \mathcal{C} \to W$ is a proper submersive holomorphic map of relative dimension 1.

$\mathcal{S}$ is a smooth complex variety, $\rho : \mathcal{S} \to W$ is a proper submersive holomorphic map of relative dimension 2.

There exists a holomorphic map $\tau : \mathcal{S} \to W \times S$, making the following diagram commutative

$$\begin{array}{ccc} \mathcal{S} & \xrightarrow{\tau} & W \times S \\ \rho \downarrow & & \mathrm{pr}_1 \downarrow \\ W & = & W \end{array}.$$

Furthermore, $\tau_{|\mathcal{S}_w} : \mathcal{S}_w \to S$ is a birational map for each $w \in W$.

Let $\phi : \mathcal{C} \to \mathcal{S}$ be a holomorphic immersion making the following diagram commutative

$$\begin{array}{ccc} \mathcal{C} & \xrightarrow{\phi} & \mathcal{S} \\ \pi \downarrow & & \rho \downarrow \\ W & = & W \end{array}.$$

Finally, let $\sigma_1, \ldots, \sigma_N$ be holomorphic sections of $\pi$, and let $m_1, \ldots, m_N$ be integers such that the zero-cycle $Z_w = \sum_i m_i \sigma_i(w)$ is of degree 0 on each component of the curve $C_w$, for each $w \in W$.

For each $i$, we get a holomorphic map

$$\alpha_i = \mathrm{pr}_2 \circ \tau \circ \phi \circ \sigma_i : W \to S,$$

and for each complex valued two-form $\omega$ on $S$, we get a two-form

$$\tilde{\omega} = \sum_i m_i \alpha_i^*(\omega)$$

462 CLAIRE VOISIN

on $W$. This two-form $\tilde{\omega}$ is Mumford's pull-back of the two-form $\omega$ on $S$ (see [7]), for the family of zero-cycles $(\text{pr}_2 \circ \tau \circ \phi(Z_w))_{w \in W}$ of $S$ parametrized by $W$.

On the other hand, for each $w \in W$, we have the Abel-Jacobi invariant $e_w := e_{C_w,S} \in J(C_w \times S)_{\text{tr}}$ or its real version

$$e_w := e_{C_w,S} \in H^1(C_w, \mathbb{Z}) \otimes_{\mathbb{Z}} H^2(S, \mathbb{Z})_{\text{tr}} \otimes_{\mathbb{Z}} \mathbb{R}/\mathbb{Z}.$$

Canonically identifying $H^1(C_w, \mathbb{Z})$ and $H^1(C_0, \mathbb{Z})$, we can view $(e_w)_{w \in W}$ as a map

$$e : W \to H^1(C_0, \mathbb{Z}) \otimes_{\mathbb{Z}} H^2(S, \mathbb{Z})_{\text{tr}} \otimes_{\mathbb{Z}} \mathbb{R}/\mathbb{Z}.$$

Clearly $e$ is differentiable (and in fact real analytic since the Abel-Jacobi invariants vary holomorphically with the parameters).

Next, for $w \in W$, the 0-cycle $Z_w$ is homologous to 0 on $C_w$, hence has a corresponding Abel-Jacobi invariant $f_w \in J(C_w)$, or its real version $f_w \in H^1(C_w, \mathbb{Z}) \otimes_{\mathbb{Z}} \mathbb{R}/\mathbb{Z}$. Identifying canonically $H^1(C_w, \mathbb{Z})$ and $H^1(C_0, \mathbb{Z})$, we can view $(f_w)_{w \in W}$ as a map

$$f : W \to H^1(C_0, \mathbb{Z}) \otimes_{\mathbb{Z}} \mathbb{R}/\mathbb{Z}.$$

Again it is easy to see that $f$ is real analytic.

Now we differentiate $e$ and $f$ to get one-forms

$$de \in \Omega_W^{\mathbb{R}} \otimes_{\mathbb{Z}} H^1(C_0, \mathbb{Z}) \otimes_{\mathbb{Z}} H^2(S, \mathbb{Z})_{\text{tr}}, \; df \in \Omega_W^{\mathbb{R}} \otimes_{\mathbb{Z}} H^1(C_0, \mathbb{Z}).$$

Finally we can contract $de \wedge df$ using the intersection pairing on $H^1(C_0, \mathbb{Z})$, to get a two-form

$$de \wedge df \in \bigwedge^2 \Omega_W^{\mathbb{R}} \otimes_{\mathbb{Z}} H^2(S, \mathbb{Z})_{\text{tr}}.$$

We can view $de \wedge df$ as an element $[de \wedge df]$ of $\text{Hom}_{\mathbb{Z}}(H^2(S, \mathbb{Z})_{\text{tr}}^{\vee}, \bigwedge^2 \Omega_W^{\mathbb{R}})$, which we can extend by $\mathbb{C}$-linearity to an element $[de \wedge df]$ of $\text{Hom}_{\mathbb{C}}(H^2(S, \mathbb{C})_{\text{tr}}^{\vee}, \bigwedge^2 \Omega_W^{\mathbb{C}})$.

Now let $\omega$ be a $(2,0)$-form on $S$, with class $[\omega] \in H^2(S, \mathbb{C})_{\text{tr}}^{\vee}$. Our main result in this section is the following:

PROPOSITION 2. *For a holomorphic two-form $\omega$ on $S$, there is the pointwise equality of two-forms on $W$*

(3.3) $$\tilde{\omega} = [de \wedge df]([\omega]).$$

The proof of formula (3.3) given below is a simplification of the original proof, following a suggestion of P. Griffiths. It goes essentially as follows: Note first that

(3.4) $$de \wedge df([\omega]) = de([\omega]) \wedge df,$$

where $de \in \text{Hom}\,(H^2(S, \mathbb{C})_{\text{tr}}^{\vee}, H^1(C_0, \mathbb{C}) \otimes \Omega_W^{\mathbb{C}})$ is the $\mathbb{C}$-linear extension of $de \in \text{Hom}\,(H^2(S, \mathbb{Z})_{\text{tr}}^{\vee}, H^1(C_0, \mathbb{Z}) \otimes \Omega_W^{\mathbb{R}})$.

Then if $\tilde{f} \in \mathcal{C}^\infty(W) \otimes H^1(C_0, \mathbb{R})$ is a lifting of $f$, we have

$$de([\omega]) \wedge df = -d(\langle de([\omega])\tilde{f}\rangle). \tag{3.5}$$

Now let $\omega'$ be the two-form on $\mathcal{C}$ induced by $\omega$ via $\mathrm{pr}_2 \circ \tau \circ \phi$. Then $\omega'$ induces a section of $\Omega_{\mathcal{C}/W} \otimes \pi^*\Omega_W$ on $\mathcal{C}$, that is a section $\beta_\omega$ of $\mathcal{H}^{1,0} \otimes \Omega_W$ on $W$. The first step is to show (see Lemma 5) that

$$de([\omega]) = \beta_\omega, \tag{3.6}$$

via the natural inclusion

$$\mathcal{H}^{1,0} \otimes \Omega_W \subset H^1_\mathbb{C} \otimes \Omega^\mathbb{C}_W \cong H^1(C_0, \mathbb{C}) \otimes \Omega^\mathbb{C}_W.$$

Next we use the definition of the Abel-Jacobi map which says that there exists a differentiably varying path $\gamma_w$ on $C_w$ such that $\partial \gamma_w = Z_w$, and for any $\eta \in H^{1,0}(C_w)$

$$\langle \eta, \tilde{f}_w \rangle = \int_{\gamma_w} \eta. \tag{3.7}$$

Combining (3.4), (3.6), and (3.7), we see that we have to show

$$\tilde{\omega} = -d(\int_\gamma \beta_\omega), \tag{3.8}$$

where $\int_\gamma \beta_\omega$ is the one-form $\psi$ on $W$ defined by $\psi(u) = \int_{\gamma_w} \beta_\omega(u)$ for $u \in T_{W,w}$. But (3.8) is essentially the homotopy formula since $\omega'$ is closed.

We now check the details of this outline of the proof and consider first the form $de$; we can view it as a map

$$[de] : H^2(S, \mathbb{Z})_{\mathrm{tr}} \to \Omega^\mathbb{R}_W \otimes_\mathbb{Z} H^1(C_0, \mathbb{Z}),$$

which can be extended by $\mathbb{C}$-linearity to a map

$$[de] : H^2(S, \mathbb{C})_{\mathrm{tr}} \to \Omega^\mathbb{C}_W \otimes_\mathbb{C} H^1(C_0, \mathbb{C}).$$

On the other hand, we have on $\mathcal{C}$ the exact sequence

$$0 \to \pi^*\Omega^2_W \to \Omega^2_\mathcal{C} \to \pi^*\Omega_W \otimes \Omega_{\mathcal{C}/W} \to 0.$$

The form $\omega' = \phi^*\tau^*\omega$ then has an image

$$\beta_\omega \in \Omega_W \otimes \mathcal{H}^{1,0}, \ \Omega_W = \Omega^{1,0}_W$$

where $\mathcal{H}^{1,0} = \pi_*\Omega_{\mathcal{C}/W}$ is the Hodge bundle with fiber $H^{1,0}(C_w) \subset H^1(C_w, \mathbb{C})$.

LEMMA 5. *For any $w \in W$, the following equality*

$$[de]([\omega])_w = (\beta_\omega)_w$$

*holds via the inclusion*

$$\Omega_{W,w} \otimes H^{1,0}(C_w) \subset \Omega_{W,w}^{\mathbb{C}} \otimes H^1(C_w, \mathbb{C}) \cong \Omega_{W,w}^{\mathbb{C}} \otimes H^1(C_0, \mathbb{C}).$$

*Proof.* Recall that $e_w \in \mathrm{Hom}\,(H^2(S, \mathbb{Z})_{\mathrm{tr}}^{\vee}, H^1(C_w, \mathbb{Z}) \otimes_{\mathbb{Z}} \mathbb{R}/\mathbb{Z})$ is obtained from the mixed Hodge structure on $H^2(S_w, C_w, \mathbb{Z})_{\mathrm{tr}}$, which fits into the exact sequence

(3.9) $\qquad 0 \to H^1(C_w, \mathbb{Z}) \to H^2(S_w, C_w, \mathbb{Z})_{\mathrm{tr}} \to H^2(S, \mathbb{Z})_{\mathrm{tr}} \to 0,$

as follows: the extension class of this extension is the class of the difference $\sigma_H - \sigma_{\mathbb{Z}} \in \mathrm{Hom}_{\mathbb{Z}}(H^2(S, \mathbb{Z})_{\mathrm{tr}}^{\vee}, H^1(C_w, \mathbb{C}))$ in the quotient

$$\mathrm{Hom}_{\mathbb{C}}(H^2(S)_{\mathrm{tr}}^{\vee}, H^1(C_w))/[F^0 \mathrm{Hom}_{\mathbb{C}}(H^2(S)_{\mathrm{tr}}^{\vee}, H^1(C_w)) \oplus \mathrm{Hom}_{\mathbb{Z}}(H^2(S, \mathbb{Z})_{\mathrm{tr}}^{\vee}, H^1(C_w, \mathbb{Z}))]$$

where $\sigma_H$ is a Hodge splitting of the sequence 3.9, and $\sigma_{\mathbb{Z}}$ is an integral splitting of the sequence 3.9. The identification

$$\mathrm{Hom}_{\mathbb{C}}(H^2(S, \mathbb{C})_{\mathrm{tr}}^{\vee}, H^1(C_w, \mathbb{C}))/F^0 \mathrm{Hom}\,(H^2(S, \mathbb{C})_{\mathrm{tr}}^{\vee}, H^1(C_w, \mathbb{C}))$$
$$\cong \mathrm{Hom}_{\mathbb{R}}(H^2(S, \mathbb{R})_{\mathrm{tr}}^{\vee}, H^1(C_w, \mathbb{R})) \cong \mathrm{Hom}_{\mathbb{Z}}(H^2(S, \mathbb{Z})_{\mathrm{tr}}^{\vee}, H^1(C_w, \mathbb{R}))$$

means simply that there is a unique splitting $\sigma_{H,\mathbb{R}}$ of the sequence 3.9 which is both Hodge and real. Then $e_w$ is the class of

$$\sigma_{H,\mathbb{R}} - \sigma_{\mathbb{Z}} \in \mathrm{Hom}_{\mathbb{Z}}(H^2(S, \mathbb{Z})_{\mathrm{tr}}^{\vee}, H^1(C_w, \mathbb{R}))$$

in the quotient $\mathrm{Hom}_{\mathbb{Z}}(H^2(S, \mathbb{Z})_{\mathrm{tr}}^{\vee}, H^1(C_w, \mathbb{Z}) \otimes_{\mathbb{Z}} \mathbb{R}/\mathbb{Z})$.

Now we have the following:

LEMMA 6. *For $\omega$ a holomorphic two-form on $S$, $\sigma_{H,\mathbb{R}}([\omega])(w)$ is the class of $\tau_w^*(\omega)$ in $H^2(S_w, C_w, \mathbb{C})_{\mathrm{tr}}$, (which is well-defined since $\tau_w^* \omega$ vanishes on $C_w$).*

*Proof.* This follows from the fact that

$$F^2 H^2(S_w, C_w)_{\mathrm{tr}} \cong F^2 H^2(S_w)_{\mathrm{tr}} \cong F^2 H^2(S)_{\mathrm{tr}},$$

so that there is a unique Hodge splitting of the sequence 3.9 over $F^2 H^2(S)_{\mathrm{tr}}$. On the other hand the map $[\omega] \mapsto$ class of $\tau_w^*(\omega)$ in $H^2(S_w, C_w, \mathbb{C})_{\mathrm{tr}}$ gives such a splitting as does $\sigma_{H,\mathbb{R}|H^{2,0}(S)}$. $\square$

Let $\mathcal{H}_{\mathbb{C}}^1$ be the flat vector bundle on $W$ with fiber $H^1(C_w, \mathbb{C})$, and $\nabla^C$ be its Gauss-Manin connection. Similarly let $\mathcal{H}_{\mathbb{C},\mathcal{S}/\mathcal{C}}^2$ be the flat vector bundle on $W$ with fiber $H^2(S_w, C_w, \mathbb{C})_{\mathrm{tr}}$, and $\nabla^{\mathcal{S}/\mathcal{C}}$ be its Gauss-Manin connection. By definition, and by Lemma 6 we have the equality:

(3.10) $\qquad [de]([\omega]) = \nabla^{\mathcal{S}/\mathcal{C}}([\tau^* \omega]),$

where $[\tau^*\omega]$ denotes the section of $\mathcal{H}^2_{\mathbb{C},\mathcal{S}/\mathcal{C}}$ whose value at $w$ is the class of $\tau^*_w\omega$ in $H^2(S_w, C_w, \mathbb{C})$. (Notice that $\nabla^{\mathcal{S}/\mathcal{C}}([\tau^*\omega])$ belongs to $\Omega^{\mathbb{C}}_W \otimes \mathcal{H}^1_{\mathbb{C}}$, since the projection of $[\tau^*\omega]$ in the quotient bundle $\mathcal{H}^2_{\mathbb{C},\mathcal{S}}$ with fiber $H^2(S_w, \mathbb{C})_{\mathrm{tr}}$ is obviously flat.) The proof of Lemma 5 follows now from the equality 3.10, and from the following general statement:

LEMMA 7. *Consider a commutative diagram of differentiable smooth fibrations*

$$\begin{array}{ccc} \mathcal{C} & \hookrightarrow & \mathcal{S} \\ \pi \downarrow & & \downarrow \rho \\ W & = & W \end{array},$$

*and let $\Omega$ be a closed $r$-form on $\mathcal{S}$, such that $\Omega_{|C_w} = 0$, for any $w \in W$. Then for the corresponding section $[\Omega]$ of the bundle $\mathcal{H}^r_{\mathcal{S}/\mathcal{C}}$, $\nabla^{\mathcal{S}/\mathcal{C}}([\Omega])$ (which belongs to $\Omega_W \otimes \mathcal{H}^{r-1}_{\mathcal{C}}/\mathcal{H}^{r-1}_{\mathcal{S}}$) can be described as follows: the restriction of $\Omega$ to $\mathcal{C}$ projects naturally to a section of $\Omega^{r-1}_{\mathcal{C}/W} \otimes \pi^*(\Omega_W)$, which is in fact vertically closed, hence gives a section $\beta_\Omega$ of $\Omega_W \otimes \mathcal{H}^{r-1}_{\mathcal{C}}$; its image in $\Omega_W \otimes \mathcal{H}^{r-1}_{\mathcal{C}}/\mathcal{H}^{r-1}_{\mathcal{S}}$ is equal to $\nabla^{\mathcal{S}/\mathcal{C}}([\Omega])$.*

*Proof.* Since the result is local, we may assume that our diagram of fibrations is trivial, that is, identifies to the inclusion $C \times W \subset S \times W$ for some $C \subset S$. For $w \in W$, $u \in T_{W,w}$, $\nabla^{\mathcal{S}/\mathcal{C}}_u([\Omega])$ is the class of the form $(d(\mathrm{int}_{\tilde{u}}\Omega) + \mathrm{int}_{\tilde{u}}(d\Omega))_{|S \times w}$, which is closed and restricts to 0 on $C_w$, in $H^r(S_w, C_w)$, where $\tilde{u}$ is the section of $T_{S \times W}$, defined along $S \times w$ and lifting $u$. Since $\Omega$ is closed, we get

$$\nabla^{\mathcal{S}/\mathcal{C}}_u([\Omega]) = \text{ class of } d(\mathrm{int}_{\tilde{u}}\Omega)_{|S \times w} \text{ in } H^r(S, C).$$

Of course $d(\mathrm{int}_{\tilde{u}}\Omega)_{|C \times w} = 0$, and the class of $\mathrm{int}_{\tilde{u}}\Omega_{|C \times w}$ in $H^{r-1}(C)$ is by definition equal to $\beta_\Omega(u)$. To conclude, it suffices to note that for an exact $r$-form $\beta = d\gamma$ on $S$ vanishing on $C$, its class in $H^{r-1}(C)/H^{r-1}(S) \subset H^r(S,C)$ is the projection of the class of $\gamma_{|C} \in H^{r-1}(C)$. So, Lemma 7, hence Lemma 5 are proved. $\square$

Now let $\tilde{f}$ be a $\mathcal{C}^\infty$ lifting of $f$ to a function with value in $H^1(C_0, \mathbb{R})$. It is clear that we have

(3.11) $$de \wedge df([\omega]) = -d(\langle de([\omega]), \tilde{f}\rangle),$$

where $\langle , \rangle$ is the intersection form on $H^1(C_0, \mathbb{C})$. Now we use the definition of the Abel-Jacobi map or Albanese map to compute this bracket; the point $\tilde{f}_w \in H^1(C_w, \mathbb{R})$ projects to

$$f^{0,1}_w \in H^{0,1}(C_w) \cong (H^{1,0}(C_w))^*$$

and we have the equality, for $\eta \in H^{1,0}(C_w))$

$$(3.12) \qquad \langle \eta, \tilde{f}_w \rangle = \langle \eta, \tilde{f}_w^{0,1} \rangle = \int_{\gamma_w} \eta,$$

for an adequate choice of a path $\gamma_w$ in $C_w$ such that $\partial \gamma_w = Z_w$.

Next, by Lemma 5, we can use this formula to compute $\langle de([\omega]), \tilde{f} \rangle$ and this gives

$$(3.13) \qquad \langle [de]([\omega]), \tilde{f} \rangle = \langle \beta_\omega, \tilde{f}^{0,1} \rangle = \int_\gamma \beta_\omega,$$

where the right-hand side is the one-form $\psi$ on $W$ defined by

$$\psi(u) = \int_{\gamma_w} \beta_\omega(u),$$

for $u \in T_{W,w}$.

By (3.11) and (3.13), to conclude the proof of Proposition 2 we have now only to prove the following:

LEMMA 8. *Let $\omega'$ be a closed holomorphic two-form on $\mathcal{C}$ with induced section $\beta_\omega$ of $\mathcal{H}^{1,0} \otimes \Omega_W$. Let $\gamma_w \subset C_w$ be a differentiably varying family of paths such that $\partial \gamma_w = Z_w$; then*

$$(3.14) \qquad \sum_i m_i \sigma_i^* \omega' = -d \left( \int_\gamma \beta_\omega \right).$$

*Proof.* It is clear that it suffices to prove equality (3.14) when we have only two sections $\sigma_1$, $\sigma_2$, and $m_1 = 1$, $m_2 = -1$. We may furthermore assume $\sigma_1(w) \neq \sigma_2(w)$ for all $w \in W$, since it suffices by continuity to prove the equality at the generic point of $W$, where this is true. (Otherwise the two sections coincide, and both sides of the equality are equal to zero.) Next, since the result is local, we can assume there is a $\mathcal{C}^\infty$ trivialization of the family $\pi: \mathcal{C} \to W$ in such a way that the two sections become constant and that there is an induced trivialization of the family of paths $\gamma_w$:

$$\begin{array}{ccc} \mathcal{C} & \cong & W \times C \\ \downarrow \pi & & \downarrow \mathrm{pr}_1 \\ W & = & W \end{array},$$

with $\sigma_i(w) = (w, c_i)$ and $\gamma_w(t) = (w, \gamma(t))$, $t \in [0,1]$, with $\gamma(0) = c_2$, $\gamma(1) = c_1$. Denote by $\Gamma: W \times [0,1] \to \mathcal{C}$ the map

$$\mathrm{Id} \times \gamma: W \times [0,1] \to W \times C \cong \mathcal{C},$$

and let $\omega'' = \Gamma^*(\omega')$. The form $\omega''$ can be written

$$(3.15) \qquad \omega'' = \eta + \delta \wedge dt,$$

where $\eta \in \mathrm{pr}_1^* \Omega_{W,\mathbb{C}}^2$, $\delta \in \mathrm{pr}_1^* \Omega_{W,\mathbb{C}}$. Then we have

$$(3.16) \qquad \sigma_1^*(\omega') - \sigma_2^*(\omega') = \eta_{|W \times 1} - \eta_{|W \times 0}.$$

Furthermore, since $\omega''$ is closed, the homotopy formula says

$$\eta_{|W \times 1} - \eta_{|W \times 0} = -d(\int_0^1 \delta_t dt). \tag{3.17}$$

Finally let $u \in T_W$ and let $\tilde{u}$ be its natural lifting in $T_{\mathcal{C}}$ given by the trivialization: then by definition of $\beta_\omega$ we have $\beta_\omega(u) = \operatorname{int}_{\tilde{u}}(\omega')_{|T_{\mathcal{C}/W}}$. Pulling this back to $W \times [0,1]$ via $\Gamma$, and using (3.15) we get

$$\delta_{w,t}(\tilde{u})dt = \gamma^*((\beta_\omega)_w(u)).$$

Fixing $w, u$ and integrating over $t$, we get

$$\int_{\gamma_w} (\beta_\omega)_w(u) = \int_0^1 \delta_{w,t}(u)dt;$$

that is,

$$\int_\gamma \beta_\omega = \int_0^1 \delta_t dt.$$

Then formula 3.14 follows from the equality above and from (3.16), (3.17). Thus Lemma 8 and Proposition 2 are proved. □

## 4. The nontriviality of $\psi_2^2$

We will prove in this section the following theorem:

THEOREM 3. *Let $S$ be a surface with $h^{2,0} \neq 0$; then $\psi_2^2(S)$ is nontrivial modulo torsion. (In fact the proof will show that $\operatorname{Im} \psi_2^2(S)$ mod. torsion is infinite dimensional.)*

The proof will be based on Propositions 3 and 4, which allow us to apply Proposition 2.

We work with the notation introduced at the beginning of Section 3, that is with the diagram

$$\begin{array}{ccccc} \mathcal{C} & \stackrel{\phi}{\to} & \mathcal{S} & \stackrel{\tau}{\to} & W \times S \\ \pi \downarrow & & \rho \downarrow & & \operatorname{pr}_1 \downarrow \\ W & = & W & = & W \end{array},$$

together with sections $\sigma_i$ of $\pi$, and integers $m_i$, defining a family of zero-cycles $Z_w$ homologous to zero on $C_w$. They allow us to define functions

$$w \mapsto e_w \in H^1(C_0, \mathbb{Z}) \otimes_\mathbb{Z} H^2(S, \mathbb{Z})_{\mathrm{tr}} \otimes_\mathbb{Z} \mathbb{R}/\mathbb{Z},$$

$$w \mapsto f_w \in H^1(C_0, \mathbb{Z}) \otimes_\mathbb{Z} \mathbb{R}/\mathbb{Z},$$

and by definition $\psi_2^2((\operatorname{pr}_2 \circ \tau \circ \phi)_* Z_w)$ is the projection modulo $U_2^2(S)$ of the product

$$e_w \cdot f_w \in H^2(S, \mathbb{Z})_{\mathrm{tr}} \otimes_\mathbb{Z} \mathbb{R}/\mathbb{Z} \otimes_\mathbb{Z} \mathbb{R}/\mathbb{Z}.$$

This product has the following explicit form: let $\{\alpha_i, \beta_i\}$, $1 \leq i \leq g$, be a symplectic basis of $H^1(C_0, \mathbb{Z})$ and let $\{\gamma_j\}$ be a basis of $H^2(S, \mathbb{Z})_{\mathrm{tr}}$; then we can write

$$e_w = \sum_{i,j} \rho_{i,j}(w) \otimes \alpha_i \otimes \gamma_j + \sum_{i,j} \chi_{i,j}(w) \otimes \beta_i \otimes \gamma_j,$$

$$f_w = \sum_i \phi_i(w) \otimes \alpha_i + \sum_i \psi_i(w) \otimes \beta_i,$$

and we have

$$(4.18) \quad e_w \cdot f_w = \sum_j (\sum_i \rho_{i,j}(w) \otimes_{\mathbb{Z}} \psi_i(w) - \sum_i \chi_{i,j}(w) \otimes_{\mathbb{Z}} \phi_i(w)) \otimes \gamma_j.$$

We prove now:

PROPOSITION 3. *Let $V \subset W$ be a smooth real analytic subset, such that for any $w \in V$, the product $e_w \cdot f_w$ vanishes in $H^2(S, \mathbb{Q})_{\mathrm{tr}} \otimes_{\mathbb{Q}} \mathbb{R}/\mathbb{Q} \otimes_{\mathbb{Q}} \mathbb{R}/\mathbb{Q}$. Then the $H^2(S, \mathbb{R})_{\mathrm{tr}}$-valued two-form $de \wedge df$ (see Section 3) vanishes on $V$.*

*Proof.* It suffices to prove that for any index $j$ the hypothesis $(e_w \cdot f_w)_j = 0$ in $\mathbb{R}/\mathbb{Q}$, for any $w \in V$, implies that $(de \wedge df)_j = 0$ on $V$. We may assume that $V$ is connected. We have then:

LEMMA 9. *There exist $I_1, I_2 \subset \{1, \ldots, g\}$, a dense subset $V' \subset V$ and coefficients*

$$\begin{aligned}
\mu_{ik} &\in \mathbb{Q}, \, i \in I_1, \, k \in \{1, \ldots, g\} - I_1, \\
\mu'_{lk} &\in \mathbb{Q}, \, l \in I_2, \, k \in \{1, \ldots, g\} - I_1, \\
\nu_{im} &\in \mathbb{Q}, \, i \in I_1, \, m \in \{1, \ldots, g\} - I_2, \\
\nu'_{lm} &\in \mathbb{Q}, \, l \in I_2, \, m \in \{1, \ldots, g\} - I_2,
\end{aligned}$$

*such that for $v \in V'$, the elements $\phi_i(v)_{i \in I_1}$, $\psi_l(v)_{l \in I_2}$ form a $\mathbb{Q}$-basis of the $\mathbb{Q}$-vector subspace of $\mathbb{R}/\mathbb{Q}$ generated by the $\phi_i(v)$, $\psi_i(v)$, $1 \leq i \leq g$. Also, the following relations hold everywhere on $V$*

$$(4.19)$$
$$\phi_k(w) = \sum_{i \in I_1} \mu_{ik} \phi_i(w) + \sum_{l \in I_2} \mu'_{lk} \psi_l(w) \text{ in } \mathbb{R}/\mathbb{Q}, \, k \in \{1, \ldots, g\} - I_1$$
$$\psi_m(w) = \sum_{i \in I_1} \nu_{im} \phi_i(w) + \sum_{l \in I_2} \nu'_{lm} \psi_l(w) \text{ in } \mathbb{R}/\mathbb{Q}, \, m \in \{1, \ldots, g\} - I_2.$$

*Proof.* Any relation $\sum_i \gamma_i \phi_i(w) + \sum_i \gamma'_i \psi_i(w) = 0$ in $\mathbb{R}/\mathbb{Q}$ holds everywhere on $V$ or only on a countable union of proper real analytic subsets of $V$. If we consider all possible such relations, it follows from Baire's theorem that there is a dense subset $V' \subset V$ such that for any $w_0 \in V'$, any relation $\sum_i \gamma_i \phi_i(w_0) + \sum_i \gamma'_i \psi_i(w_0) = 0$ in $\mathbb{R}/\mathbb{Q}$ implies that $\sum_i \gamma_i \phi_i(w) + \sum_i \gamma'_i \psi_i(w) = 0$ in $\mathbb{R}/\mathbb{Q}$, for any $w \in V$. Choosing for such $w_0$ a basis $\phi_i(w_0)$, $\psi_l(w_0)$, $i \in I_1$, $l \in I_2$ of the $\mathbb{Q}$-vector subspace of $\mathbb{R}/\mathbb{Q}$ generated by the $\phi_i(w_0)$, $\psi_i(w_0)$, $1 \leq i \leq g$, gives the result. $\square$

Now in formula (4.18) we replace $\phi_k(w)$ and $\psi_m(w)$ by their expressions in (4.19), which gives

$$
\begin{aligned}
(e_w \cdot f_w)_j &= \sum_{l \in I_2} \rho_{lj}(w) \otimes \psi_l(w) + \sum_{m \notin I_2} \rho_{mj}(w) \\
&\otimes (\sum_{i \in I_1} \nu_{im}\phi_i(w) + \sum_{l \in I_2} \nu'_{lm}\psi_l(w)) \\
&- \sum_{i \in I_1} \chi_{ij}(w) \otimes \phi_i(w) - \sum_{k \notin I_1} \chi_{kj}(w) \\
&\otimes \sum_{i \in I_1} \mu_{ik}\phi_i(w) + \sum_{l \in I_2} \mu'_{lk}\psi_l(w)),
\end{aligned}
$$

where the equality holds in $\mathbb{R}/\mathbb{Q} \otimes_{\mathbb{Q}} \mathbb{R}/\mathbb{Q}$.

Now we use the fact that $\phi_i(w)$ and $\psi_l(w)$ are independent over $\mathbb{Q}$ for $w$ in the dense subset $V'$. Then for $w \in V'$ the condition $(e_w \cdot f_w)_j = 0$ in $\mathbb{R}/\mathbb{Q} \otimes_{\mathbb{Q}} \mathbb{R}/\mathbb{Q}$ implies

$$
(4.20) \quad \rho_{lj}(w) + \sum_{m \notin I_2} \nu'_{lm}\rho_{mj}(w) - \sum_{k \notin I_1} \mu'_{lk}\chi_{kj}(w) = 0 \quad \text{in } \mathbb{R}/\mathbb{Q}, \forall l \in I_2
$$

$$
-\chi_{ij}(w) + \sum_{m \notin I_2} \nu_{im}\rho_{mj}(w) - \sum_{k \notin I_1} \mu_{ik}\chi_{kj}(w) = 0 \quad \text{in } \mathbb{R}/\mathbb{Q}, \forall i \in I_1.
$$

But recall that $V'$ is the complementary set in $V$ of a countable union of proper real analytic subsets. So the equalities (4.20), being satisfied on $V'$, must hold everywhere on $V$.

Now we can differentiate (4.19) and (4.20): indeed these equalities mean that for liftings of the functions $\phi_i, \psi_i, \chi_{ij}, \rho_{kj}$ to functions with values in $\mathbb{R}$, the corresponding equalities hold modulo some (necessarily constant) rational numbers. This gives

$$
(4.21) \quad d\phi_k = \sum_{i \in I_1} \mu_{ik}d\phi_i + \sum_{l \in I_2} \mu'_{lk}d\psi_l, \, k \in \{1,\ldots,g\} - I_1
$$

$$
d\psi_m = \sum_{i \in I_1} \nu_{im}d\phi_i + \sum_{l \in I_2} \nu'_{lm}d\psi_l, \, m \in \{1,\ldots,g\} - I_2.
$$

$$
(4.22) \quad d\rho_{lj} + \sum_{m \notin I_2} \nu'_{lm}d\rho_{mj} - \sum_{k \notin I_1} \mu'_{lk}d\chi_{kj} = 0, \quad \text{for all } l \in I_2
$$

$$
-d\chi_{ij} + \sum_{m \notin I_2} \nu_{im}d\rho_{mj} - \sum_{k \notin I_1} \mu_{ik}d\chi_{kj} = 0, \quad \text{for all } i \in I_1.
$$

Next we have

$$(de \wedge df)_j = (de)_j \wedge df = \sum_i d\rho_{ij} \wedge d\psi_i - \sum_i d\chi_{ij} \wedge d\phi_i,$$

which shows that this is equal to zero, by (4.21) and (4.22). Now Proposition 3 is proved. □

Combining Proposition 2 and Proposition 3, we conclude:

COROLLARY 1. *Under the assumptions of Proposition 3, the pull-back*

$$\tilde{\omega} = \sum_i m_i (\mathrm{pr}_2 \circ \tau \circ \phi \circ \sigma_i)^*(\omega)$$

*of any holomorphic two-form $\omega$ on $S$ vanishes on $V$.*

Next we have the following:

PROPOSITION 4. *Assume the map $\psi_2^2(S)$ vanishes modulo torsion in $J_2^2(S)$; then there exist data*

$$\begin{array}{ccccc} \mathcal{C} & \stackrel{\phi}{\to} & \mathcal{S} & \stackrel{\tau}{\to} & W \times S \\ \pi \downarrow & & \rho \downarrow & & \mathrm{pr}_1 \downarrow \\ W & = & W & = & W \end{array},$$

*together with sections $\sigma_i$ of $\pi$, and integers $m_i$, defining a family of zero-cycles $Z_w$ homologous to zero on $C_w$ satisfying the following properties:*

- *There exists a map $\psi = (\psi_1, \psi_2) : W \to S \times S$ such that*

$$(\mathrm{pr}_2 \circ \tau \circ \phi)_* Z_w = \psi_1(w) - \psi_2(w)$$

  *as a zero-cycle of $S$, for any $w \in W$.*

- *There is a smooth locally closed real analytic subset $V \subset W$ such that for any $w \in V$, $e_w \cdot f_w$ vanishes in $H^2(S, \mathbb{Q})_{\mathrm{tr}} \otimes_\mathbb{Q} \mathbb{R}/\mathbb{Q} \otimes_\mathbb{Q} \mathbb{R}/\mathbb{Q}$, and $\psi_{|V}$ is a submersion.*

Clearly this proposition implies Theorem 3; indeed, if $\psi_2^2(S)$ vanishes modulo torsion in $J_2^2(S)$, Proposition 4 and Corollary 1 give a submersive map $\psi : V \to S \times S$ such that for any holomorphic two-form $\omega$ on $S$,

$$\psi_1^*(\omega) - \psi_2^*(\omega) = \sum_i m_i (\mathrm{pr}_2 \circ \tau \circ \phi \circ \sigma_i)^*(\omega)$$

vanishes on $V$. It follows that $\mathrm{pr}_1^*(\omega) - \mathrm{pr}_2^*(\omega)$ vanishes on an open set of $S \times S$, hence that $\omega = 0$. So we have proved that $\psi_2^2(S) = 0$ modulo torsion implies $H^{2,0}(S) = \{0\}$, that is, Theorem 3. □

*Proof of Proposition* 4. By definition of $\psi_2^2$, the assumption implies that for any $(x_1, x_2) \in S \times S$, there exist a smooth curve $C$, a zero-cycle $Z$ homologous to zero on $C$, and an immersion $\phi : C \hookrightarrow \tilde{S}$ of $C$ in a surface $\tilde{S} \xrightarrow{\tau} S$ birational to $S$, such that $\tau \circ \phi_*(Z) = x_1 - x_2$ and $e_{C,S} \cdot f_{Z,C} = 0$ in $H^2(S, \mathbb{Q})_{\text{tr}} \otimes_\mathbb{Q} \mathbb{R}/\mathbb{Q} \otimes_\mathbb{Q} \mathbb{R}/\mathbb{Q}$. Now note that there are countably many quasi-projective varieties (that we may assume smooth by desingularization) $W_m$, together with data

$$\begin{array}{ccccc} \mathcal{C}_m & \xrightarrow{\phi} & \mathcal{S}_m & \xrightarrow{\tau} & W_m \times S \\ \pi_m \downarrow & & \rho_m \downarrow & & \text{pr}_1 \downarrow \\ W_m & = & W_m & = & W_m \end{array},$$

with sections $\sigma_i^m$ of $\pi_m$, and integers $m_i^m$, defining a family of zero-cycles $Z_w^m$ homologous to zero on $C_w^m$, and satisfying the following properties:

- There exists a map $\psi_m = (\psi_1^m, \psi_2^m) : W_m \to S \times S$ such that $(\tau_m \circ \phi_m)_* Z_w^m = \psi_1^m(w) - \psi_2^m(w)$ as a zero-cycle of $S$, for any $w \in W_m$.

- Any set of data $((x_1, x_2), C, Z, \phi, \tau)$ as above, such that $\tau \circ \phi_*(Z) = x_1 - x_2$ identifies with the data parametrized by some point $w \in W_m$, with $(x_1, x_2) = \psi_m(w)$.

On each $W_m$, we have the locally defined maps

$$\begin{array}{rcl} e_m : W_m & \to & H^1(C_0^m, \mathbb{Z}) \otimes_\mathbb{Z} H^2(S, \mathbb{Z})_{\text{tr}} \otimes_\mathbb{Z} \mathbb{R}/\mathbb{Z}, \\ f_m : W_m & \to & H^1(C_0^m, \mathbb{Z}) \otimes_\mathbb{Z} \mathbb{R}/\mathbb{Z}, \end{array}$$

(which are globally defined as sections of a flat bundle), and their product

$$e_m \cdot f_m : W_m \to H^2(S, \mathbb{Q})_{\text{tr}} \otimes_\mathbb{Q} \mathbb{R}/\mathbb{Q} \otimes_\mathbb{Q} \mathbb{R}/\mathbb{Q}.$$

So the assumption of Proposition 4 is that

$$S \times S = \bigcup_m \psi_m((e_m \cdot f_m)^{-1}(0)).$$

We have now:

LEMMA 10. *Locally $(e_m \cdot f_m)^{-1}(0))$ is a countable union of real analytic subsets of $W_m$.*

Assuming Lemma 10 we have countably many locally closed real analytic subsets $W_m^n \subset W_m$ on which $e_m \cdot f_m$ vanishes, and such that

$$S \times S = \bigcup_{m,n} \psi_m(W_m^n).$$

Stratifying each $W_m^n$ into smooth real analytic subsets, we may assume the $W_m^n$ are smooth. The theorems of Baire and Sard imply now that for some $(m, n)$, $\psi_{m|W_m^n}$ must be submersive at some point of $W_m^n$, hence on an open subset $V$ of it. So Proposition 4 is proved, with $W = W_m$, $\psi = \psi_m$. □

*Proof of Lemma* 10. The proof was almost completed in the course of the proof of Proposition 3. With the notation introduced there (and forgetting the subscript $m$), it follows from the computations made there that, for the $j^{\text{th}}$-component $(e \cdot f)_j$ of $e \cdot f$, $(e \cdot f)_j^{-1}(0) \subset W$ can be written locally as the countable union of the sets $W_{I_1, I_2, \mu_{ik}, \mu'_{lk}, \nu_{im}, \nu'_{lm}} \subset W$ where the equations (4.19) and (4.20) are satisfied. But choosing (locally) liftings of the $\phi_i, \psi_i, \rho_{ij}, \chi_{ij}$ to real analytic functions with values in $\mathbb{R}$, one sees immediately that each $W_{I_1, I_2, \mu_{ik}, \mu'_{lk}, \nu_{im}, \nu'_{lm}}$ is a countable union of real analytic subsets of $W$. Now the lemma is proved, since $(e \cdot f)^{-1}(0) = \bigcap_j (e \cdot f)_j^{-1}(0)$. □

*Remark* 1. More generally, we have proved that in $S^{[k]} \times S^{[k]}$, the set $Z$ of points $(z_1, z_2)$ such that $\psi_2^2(z_1 - z_2) = 0$ mod torsion is covered by a countable union of images of real analytic sets $V, \psi : V \to S^{[k]} \times S^{[k]}$, such that for any holomorphic two-form $\omega$ on $S$ with induced form $\omega_k$ on $S^{[k]}$, $\psi_1^* \omega_k - \psi_2^* \omega_k$ vanishes on $V$. Hence, the Mumford argument (see [7]) applies to show that $\operatorname{Im} \psi_2^2$ mod torsion is infinite dimensional. Indeed, if $z \in S^{[k]}$ is a general point, the two-form $\psi_1^* \omega_k - \psi_2^* \omega_k$ is nondegenerate at $(z, z)$, so that its real part is also nondegenerate, and the fact that it vanishes on $Z$ implies that the real dimension of any component of $Z$ passing through $(z, z)$ is at most equal to $\frac{\dim_{\mathbb{R}} S^{[k]} \times S^{[k]}}{2} = \dim_{\mathbb{R}} S^{[k]}$. But any component of $Z$ passing through $(z, z)$ has to dominate an open set of $S^{[k]}$ by the first projection, if $z$ is chosen outside a countable union of real analytic sets. It follows that the map $z' \mapsto \psi_2^2(z' - z)$ has almost all of its fibers countable in some neighbourhood of $z$. Hence the dimension of its image, defined as $\dim_{\mathbb{R}} S^{[k]} - \dim_{\mathbb{R}} (\textit{general fiber})$, is equal to $\dim_{\mathbb{R}} S^{[k]}$, which tends to $\infty$ with $k$. Here *general* is with respect to the real analytic Zariski topology and the dimension of a fiber is well-defined since the fiber is covered by a countable union of real analytic sets.

Université Pierre et Marie Curie, Paris, France
*E-mail address*: voisin@math.jussieu.fr